\documentclass[final,twoside, a4paper]{amsart}

\usepackage[english]{babel}
\usepackage{amsmath, amsfonts,amssymb,amsopn,amscd,amsthm}
\usepackage{mathpazo}
\usepackage{graphicx}

 \usepackage{hyperref}
 \hypersetup{colorlinks=true,linkcolor=blue,citecolor=blue,linktoc=page}
 
\numberwithin{equation}{section}

\renewcommand{\ni}{\noindent}
\newcommand{\eq}{\begin{equation}}
\newcommand{\qe}{\end{equation}}
\renewcommand{\d}{\mathrm{d}}
\newcommand{\T}{\mathbb{T}}
\renewcommand{\d}{\mathrm{d}}
\newcommand{\R}{\mathbb{R}}                     
\newcommand{\M}{\mathcal{M}}
\def\P{\mathbb{P}}
\def\E{\mathbb{E}} 
\def\Var{\mathop{\rm Var}}
\newcommand{\normTV}[1]{\left\lVert #1\right\lVert_{TV}}
\newcommand{\abs}[1]{\left\vert#1\right\vert}

\theoremstyle{plain}
\newtheorem{thm}{Theorem}[section]
\newtheorem{lem}[thm]{Lemma}
\newtheorem{prop}[thm]{Proposition}
\newtheorem{cor}{Corollary}
\newtheorem{defn}{Definition}[section]
\theoremstyle{definition}

\theoremstyle{remark}
\newtheorem*{rem}{Remark}

\begin{document}
\sloppy

\pagestyle{headings} 
\title{Spike detection from inaccurate samplings}
\date{\today}

\author{Jean-Marc Aza\"is \and Yohann de Castro \and Fabrice Gamboa}
\address{JMA and FG are with the Institut de Math\'ematiques de Toulouse (CNRS UMR 5219). Universit\'e Paul Sabatier, 118 route de Narbonne, 31062 Toulouse, France.}
\email{$\{$jean-marc.azais$\}$,$\{$fabrice.gamboa$\}$@math.univ-toulouse.fr}
\urladdr{www.math.univ-toulouse.fr/~$\{$azais$\}$,$\{$gamboa$\}$}
\address{YdC is with the D\'epartement de Math\'ematiques (CNRS UMR 8628), B\^atiment 425, Facult\'e des Sciences d'Orsay, Universit\'e Paris-Sud 11, F-91405 Orsay Cedex, France.}
\email{yohann.decastro@math.u-psud.fr}
\urladdr{www.math.u-psud.fr/$\sim$decastro} 

\begin{abstract}
This article investigates the support detection problem using the LASSO estimator in the space of measures. More precisely, we study the recovery of a discrete measure (spike train) from few noisy observations (Fourier samples, moments...) using an $\ell_{1}$-regularization procedure. In particular, we provide an explicit quantitative localization of the spikes. 
\end{abstract}

\maketitle 

\section{Introduction}

\subsection{Super-resolution}

Imaging experiments can be subject to device limitations where one cannot observe enough information in order to recover fine details. For instance, in optical imaging, the physical limitations are evaluated by the resolution. This latter measures the minimal distance between lines that can be distinguished. Hence, the details below the resolution limit seem unreachable. The super-resolution phenomenon is the ability to recover the information beyond the physical limitations. Surprisingly, if the object of interest is simple, e.g. a discrete measure, then it is possible to override the resolution limit. In particular, the reader may think of important questions in applied harmonic analysis such as the problem of source separation. Many companion applications in astronomy, medical imaging and single molecule imaging in 3D microscopy are at stake, see \cite{park2003super, lindberg2012mathematical,studer2012compressive} and references therein. Hence, theoretical guarantees of source detection are of crucial importance in practice. 

In this paper, we prove quantitative detection guarantees from noisy observations (Fourier samples, moments samples). Furthermore, these quantitative estimates can be computed using a tractable algorithm, called BLASSO.

\subsection{Previous works}

The theoretical analysis of the $\ell_{1}$-regularization in the space of measures was initiated by Donoho \cite{donoho1992superresolution}. Few years after, rates of convergence in super-resolution have been investigated by P. Doukhan, E. Gassiat and one author of this present paper \cite{MR1393035,MR1393430}. They considered the exact reconstruction of a nonnegative measure and derived results when one only knows the values of a finite number of linear functionals at the target measure. Moreover, they study stability with respect to a metric for weak convergence. Likewise, two authors of this paper \cite{de2012exact} proved that $k$ spikes trains can be faithfully resolved from $m=2k+1$ samples (Fourier samples, Stieltjes transformation, Laplace transform...) by using an $\ell_{1}$-minimization method.

Following recent proposal \cite{bredies2010inverse, burger2004convergence,resmerita2006error,scherzer2009variational} on inverse problems regularization in Banach spaces, we consider convergence rates in Bregman divergence. On a more general note, inverse problems on the space of measures are now well understood, see \cite{hofmann2007convergence,scherzer2009sparsity} for instance. We capitalize on these earlier works to construct our analysis. In particular, we use them to give quantitative localizations of the recovered spikes, which is new.

In the super-resolution frame, the important paper \cite{candes2012towards} shows that if the spikes are well ``seperated'' then there exists a dual certificate, i.e. an $\ell_{\infty}$-constrained trigonometric polynomial that interpolates the phase of the weights at the spikes locations. This construction provides a quadratic isolation condition, see Definition \ref{def:QIC}, of the spikes which is crucial for $\ell_{1}$-minimization in the space of measures. In a predating paper \cite{de2012exact}, the authors investigates $\ell_{1}$-minimization with different types of measurements: trigonometric, polynomial, Laplace transform... In view of application, the recent works \cite{candes2012super,tang2013near} derive results in $\ell_{1}$ and $\ell_{2}$ prediction, i.e. the estimation of the input frequencies. Moreover, note that noise robustness of support recovery is proved in \cite{duval2013exact}.

A postdating paper \cite{fernandez2013support} bounds the support detection error for a constrained formulation of the $\ell_{1}$-minimization in the space of measures as in Theorem \ref{thm:LocBIPGeneral}. Regarding unconstrained formulation, a second postdating paper \cite{tang2013near} studies spikes detection for the Fourier sampling under a Gaussian noise model. The authors provide an optimal rate for the $\ell_{2}$-prediction, namely the $\ell_{2}$ distance between the recovered Fourier coefficients and the original ones. 

The aforementioned results suggest that the recovered spike locations should be close to the input support. This is investigated, for the first time, in this paper for an unconstrained $\ell_{1}$-minimization problem under a general sampling scheme. 

\subsection{General model and notation}

Let $({\mathbb T},d)$ be a compact metric space homeomorphic to either the interval $([0,1],\lvert\,.\,\lvert)$ or the unit circle $\mathbb S^1$, which is identified to the metric space $(\R\mod(1),d(\cdot,\cdot))$ via the mapping $z = e^{\mathrm i2\pi t}$. In this latter case, the distance $d$ is taken around the circle. Let $\Delta$ be a complex measure on ${\mathbb T}$ with discrete support of size $s$. In particular, the measure $\Delta$ has polar decomposition, see \cite{rudin1987real} for a definition:
\begin{equation}\label{def:TargetMeasure}
\Delta=\sum_{k=1}^s\Delta_k\,\exp(\mathrm{i}\,\theta_k)\,\delta_{T_k}\,,
\end{equation}
 where $\Delta_k>0$, $\theta_k\in\R$, $T_k\in {\mathbb T}$ for $k=1,\ldots,s$ and
$\delta_x$ denotes the Dirac measure at point $x$. Let $m$ be a positive integer and $\mathbf F=\{\varphi_0,\varphi_1,\dotsc,\varphi_m\}$ be a family of complex continuous functions on ${\mathbb T}$. Assume that the family $\mathbf F$ is \textbf{a family of orthonormal functions} with respect to a probability measure $\Pi$ on $\mathbb T$. 

We recall some basic concepts in the frame of generalized polynomials. Define the $k$-th {generalized moment} of a complex measure $\mu$ on ${\mathbb T}$ as:
\[
c_{k}(\mu)=\displaystyle\int\nolimits_{\mathbb T} \varphi_{k}\,\d\mu\,, 
\]
for $k=0,1,\dotsc,m$. Suppose we observe $y=(y_k)_{k=0}^m$ defined by $y_k=c_k(\Delta)+\varepsilon_k$ for $k=0,1,\dotsc,m$, where $\varepsilon=(\varepsilon_k)_{k=0}^m$ is a complex valued white noise. This can be written as:
\[
 y=\int_{\mathbb T}\Phi\,\d\Delta+\varepsilon\,,
\]
where $\Phi=(\varphi_0,\ldots,\varphi_m)$. In this paper, we investigate the recovery of the complex measure $\Delta$ from $m+1$ measurements given by $y$. 

\begin{rem}
 Along this article, we shall mention examples in the Fourier case (Fourier samples) or in the polynomial case (moment samples), notation are described therein. If not specified, notation are in accordance with the general model.
\end{rem}

\subsection{Unconstrained minimization}

Denote by $\M$ the set of all finite complex measures on ${\mathbb T}$ and by $\normTV{\,.\,}$ the {total variation} norm. We recall that for all $\mu\in\M$, \[\normTV{\mu}=\sup_{\mathbf P}\sum_{E\in\mathbf P}\abs{\mu(E)}\,,\] where the supremum is taken over all partitions $\mathbf P$ of ${\mathbb T}$ into a finite number of disjoint measurable subsets. For further details, we refer the reader to \cite{rudin1987real}. 
%

\ni
By analogy with the LASSO \cite{MR1379242}, \textit{Beurling LASSO} (BLASSO) is the process of reconstructing a discrete measure $\Delta$ from the samples $y$ by finding a solution to:

\eq
\label{def:Blasso}
\tag{BLASSO}
\hat\Delta\in\arg\min_{\mu\in\M}\frac12\lVert
\int_{\mathbb T}\Phi\,\d\mu-y\lVert^2_2+\lambda\lVert\mu\lVert_{TV}\,,
\qe

\ni
where $\lambda$ is a tuning parameter. 

\begin{rem}
 For the case of Fourier coefficients and $\varepsilon=0$, \eqref{def:Blasso} is called {\it Generalized Minimal Extrapolation} (GME), see \cite{de2012exact} and Eq. \eqref{support pursuit} for a definition. This procedure finds an extrapolation of a function given on a subset of $\mathbb T$ with minimal TV-norm of its Fourier transform among all possible extrapolating functions. Moreover, our procedure is an extension of the LASSO estimator of the high-dimensional regression theory.
\end{rem}

\begin{rem}
On the algorithmic side, observe that the Fenchel dual program of \eqref{def:Blasso} can be recast into a SDP program when using Fourier samplings or moment samplings, see \cite{candes2012towards,tang2013near} for instance.
\end{rem}

\ni
One knows that the extreme points of the unit ball of the TV-norm are the atoms $\delta_x$ where $x\in\mathbb T$. Therefore, $\ell_{1}$-minimization compels the solutions to be discrete measures. Nevertheless, the TV-norm is innappropriated in measuring the distance between spikes, and it seems difficult to localize BLASSO.  Two questions immediately arise:

\begin{itemize}
 \item How close is the recovered spike train $\hat\Delta$ from the target $\Delta$?
  \item How accurate is the localization of BLASSO in terms of the noise and the amplitude of the recovered/original spike?
\end{itemize}

\ni
To the best of our knowledge, this paper is the first to address these questions in a general frame. In particular, it is the first paper to provide quantitative localization guarantees from the output amplitudes.

\subsection{Contribution} 

The present paper is concerned with the problem of the localization of the target spikes. In the general frame, we show that the large recovered spikes are very close to the target spikes, and that the mass of all the recovered spikes far from the original support is small. In the Fourier case and the moment case, we specify that under mild assumptions on the support of the target, the mass of the reconstructed measure is concentrated around large spikes, i.e. spikes above the noise level. In particular, we derive error bounds in terms of the recovered measure and the original measure.

\subsection{Organization of the paper}

The next section present the main results in the general frame. We derive the corresponding results in the Fourier frame and the moment frame in Section 3. The bound on the noise in the space of continuous functions is given in Section 4 using the Rice method. Section 5 recalls some useful optimization tools in order to implement our procedure.

\section{Quantitative localization in the general frame}

\subsection{Zero-noise problem}

All solution to $\ell_{1}$-regularization satisfies some optimality condition based on the fact that the sub-gradient of the regularization function vanishes at the solution point. Then a sufficient condition for exact recovery is that the target measure satisfies this optimality condition. This analysis has led to the notion of dual certificate, see \cite{de2012exact,candes2012towards} for instance, and the notion of source condition, see \cite{burger2004convergence}.

\begin{defn}[Dual certificate]
We say that a generalized polynomial $P=\sum_{k=0}^{m}a_k\varphi_{k}$ is a dual certificate for the measure $\Delta$ defined by \eqref{def:TargetMeasure} if and only if it satisfies the following properties:
 \begin{itemize}
  \item phase interpolation:  $\forall k\in\{1,\ldots,s\}\,,\quad P(T_k)=\exp(-\mathbf{i}\theta_k)$,
  \item $\ell_{\infty}$-constraint: $\lVert P\lVert_{\infty}\leq1$.
 \end{itemize}
\end{defn}

\ni
Indeed, one can prove that the target measure $\Delta$ is a solution of the following $\ell_{1}$-minimization:
\begin{equation}\label{support pursuit}
 \Delta^{GME}\in\arg\min\displaylimits_{\mu\in\mathcal M}\lVert{\mu}\lVert_{TV}\quad \text{s.t.} \
\int_{\mathbb T} \Phi\,\d\mu=\int_{\mathbb T} \Phi\,\d\Delta\,.
\end{equation}
\ni
{if and only if} $\Delta$ has a dual certificate, a proof can be found in \cite{de2012exact}. We understand that if $\Delta$ has no dual certificate (i.e. no generalized polynomial with infinity norm less than $1$ interpolates the phases at support point locations) there is no hope in recovering $\Delta$ with an $\ell_{1}$-minimization method.

\subsection{Quantitative localization guarantees from the output amplitudes}\label{subsec:separable}
In the presence of noise, one cannot ask for exact recovery and the notion of dual certificate is too loose for establishing stability results. Therefore, we strengthen a bit more this notion so as to derive a quantitative localization of the target spikes from BLASSO. 

 \begin{defn}[Quadratic isolation condition] 
 \label{def:QIC}
 We say that a finite set $\mathbf S=\{T_1,\ldots,T_s\}\subset\mathbb T$ satisfies the quadratic isolation condition with parameters $C_{a}>0$ and $0<C_{b}<1$, denoted by $\mathrm{QIC}(C_{a}, C_{b})$, if and only if for all $(\theta_k)_{k=1}^{s}\in\mathbb R^{s}$, there exists $P\in\mathrm{Span}(\mathbf F)$ such that  for all $k=1,\ldots,s$, $P(T_k)=\exp(-\mathbf{i}\theta_k)$, and
\[\forall x\in\T\,,\quad 1-\lvert P(x)\lvert\geq \min_{T\in\mathbf S}\{C_{a}m^{2}d(x,T)^{2},C_{b}\}\,.\]
\end{defn}

\ni
Before stating the theorem, we would like to point out that, under general conditions \cite{krein1977markov}, the solutions of the convex program \eqref{def:Blasso} always contain an atomic solution with support of size less than $m+2$. Moreover, this solution can be computed from the solution of the convex dual problem in the Fourier frame, see \cite{candes2012towards} for instance. Therefore, we can always consider an atomic solution of BLASSO:

\[
\hat\Delta=\sum_{k=1}^{n}\hat\Delta_k\,\exp(\mathrm{i}\,\hat\theta_k)\,\delta_{\hat T_k}\,.
\]

\ni 
We begin with a first result showing that, even with a small regularization parameter $\lambda$, the large spikes of BLASSO detect the original support points. Thereafter, we will denote
\[
 \mathcal N(c_{0},\mathbf S):=\Big\{x\in\mathbb T\,;\quad \min_{T\in\mathbf S}d(x,T)\leq\frac{c_0}{m}\Big\}\,,
 \]
 and $\mathcal F(c_{0},\mathbf S)$ its complement, or $ \mathcal N$ and $\mathcal F$ for short.

\begin{thm}
\label{thm:LocQICGeneral}
 Assume that the support $\mathbf S$ of $\Delta$ satisfies $\mathrm{QIC}(C_{a}, C_{b})$. Let $\lambda$ be such that $\lambda\geq\lVert\varepsilon\lVert_{2}$. Then, for all $x\in\mathbb T$ such that: 
\[
 \lvert\hat\Delta(\{x\})\lvert>\dfrac{2\lambda}{C_b}\,,
\]
there exists $T\in\mathbf S$ satisfying:
\[
 d(x,T)\leq \Big[\dfrac{2\lambda}{C_a\lvert\hat\Delta(\{x\})\lvert}\Big]^{1/2}\,\dfrac1{m}\leq\dfrac{c_{0}}{m}\,,
\]
where $c_{0}=\sqrt{C_{b}/C_{a}}$. Observe that, if $\displaystyle\min_{k\neq l}d(T_{k}-T_{l})>\frac{2c_{0}}m$ then the aforementioned point $T$ is unique. Moreover,
\begin{enumerate}
\item
$\displaystyle\sum_{\hat T_{k}\in \mathcal N(c_{0},\mathbf S)}\hat\Delta_{k}\,\min_{T_{l}\in\mathbf S}d(\hat T_{k},T_l)^{2}\leq \frac{2\lambda}{C_{a}m^{2}}$,
\item
$\displaystyle\sum_{\hat T_{k}\in \mathcal F(c_{0},\mathbf S)}\hat\Delta_{k}\leq \frac{2\lambda}{C_{b}}$.
\end{enumerate}

\end{thm}

\ni
This result shows that BLASSO puts a small weight far from the original support, see $(2)$. Moreover, it shows that the reconstructed points, with large weights, are close to the true support, see $(1)$. Our result quantitatively bounds the support recovery error in terms of the amplitude of the solution. In actual practice, the result $(1)$ provides a confidence set on the localization of the true support. Furthermore, it shows that this localization is getting better as the recovered amplitude is large.

\begin{rem} Postdating this paper, some important works \cite{fernandez2013support,tang2013near} improve our result, in the Fourier case, providing bounds that depend only on the amplitude of the original spike. An attentive reader can see that the bounds $(ii)$ and $(iii)$ in Theorem 1.2 of \cite{fernandez2013support} and the bounds $(i)$ and $(ii)$ in Theorem 2 of \cite{tang2013near} are covered by our predating result. Moreover, Theorem \ref{thm:LocQICGeneral} deals with a regularizing parameter $\lambda$ that can be of the order of $\lVert\varepsilon\lVert_{2}$. In contrast, the regularizing parameter $\tau$ in \cite{tang2013near} has to be of the order of $\lambda_{0}=\lVert \sum_{k=0}^{m}\varepsilon_k\varphi_{k}\lVert_{\infty}$ which is roughly of the order of $\lVert \varepsilon\lVert_{\infty}$ in the Gaussian noise model, see \ref{prop:RiceFourier}.
\end{rem}

\subsection{Quantitative localization guarantees from the input amplitudes} 

One knows \cite{borwein1995polynomials} that any bounded polynomial on a compact of the real line has its derivative upper bounded by a constant times its degree. Hence, a uniform upper bound on the derivatives of the dual certificate $P$ can be given. 
\begin{defn}[Bernstein Isolation Property]
We say that a set $\mathbf S$ satisfies the Bernstein Isolation Property with parameters $c_0>0$ and $C_{c}>0$, denoted by $\mathrm{BIP}(c_0, C_{c})$, if and only if 
\[
\forall P\in\mathrm{Span}(\mathbf F),\quad\forall x\in\mathcal N(c_{0},\mathbf S)\,, \quad\lvert P''(x)\lvert\leq C_{c}m^{2}\lVert P\lVert_{\infty}\,.
\]
 \end{defn} 
\ni
Using this property, we can specify the result of Theorem \ref{thm:LocQICGeneral} in terms of the input amplitudes.
\begin{thm}
\label{thm:LocBIPGeneral}
Assume that the support $\mathbf S$ of $\Delta$ satisfies $\mathrm{BIP}(c_0, C_{c})$ and $\mathrm{QIC}(C_{a}, C_{b})$ with $c_{0}=\sqrt{C_{b}/C_{a}}$. Let $\lambda$ be greater than $\lVert \sum_{k=0}^{m}\overline\varepsilon_k\varphi_{k}\lVert_{\infty}$. Then
\begin{enumerate}
\setcounter{enumi}{2}
\item 
$\displaystyle\lvert \Delta_{j}-\sum_{\{k:\ d(\hat T_{k},T_{j})\leq\frac{c_{0}}m\}} \hat\Delta_{k}\lvert\leq C'\lambda\,,$
\end{enumerate}
where $C'=2+\max\Big\{\frac{2(1-C_{b})}{C_{b}},\frac{C_{c}}{C_{a}}\Big\}$. Moreover for all $T_{j}\in \mathbf S$ corresponding to a weight such that: 
\[
 \Delta_{j}>C'\lambda\,,
\]
there exists $\hat T\in\mathrm{Supp}(\hat\Delta)$ satisfying:
\[
 d(T_{j},\hat T)\leq \Big[\dfrac{2\lambda}{C_a(\Delta_{j}-C'\lambda)}\Big]^{1/2}\,\dfrac1{m}\,.
\]
\end{thm}

\ni 
Observe that this phenomenon has been investigated in the Fourier case, see for instance \cite{fernandez2013support,tang2013near}. Theorem \ref{thm:LocBIPGeneral} shows that it is conducted by the BIP property, namely a control of the second derivative of the dual certificate on a small vicinity of each point of the support $\mathbf S$. Hence, it may applies to a more general frame.

\subsection{Examples of families satisfying BIP}

The BIP property is satisfied for a large class of family of measurements, see \cite{borwein1995polynomials} for instance. We present some standard Bernstein-type inequalities for which the BIP property can be derived.

\begin{description}
\item[Fourier samples] Let $\mathbf F_F=\left\{\exp(-\imath2\pi f_{c} x),\dotsc,1,\dotsc,\exp( \imath 2\pi f_{c} x)\right\}$  then 
\eq\label{eq:BernsteinFourier}
\forall P\in\mathrm{Span}(\mathbf F_{F}),\quad\forall x\in[0,1]\,,\quad\lvert P''(x)\lvert\leq \pi^{2}(2f_{c})^{2}\lVert P\lVert_{\infty}\,.
\qe
\item[Moment samples] Let $\mathbf F_M=\left\{1,x,\dotsc,x^{m}\right\}$  then 
\eq\label{eq:BernsteinMoment}
\forall P\in\mathrm{Span}(\mathbf F_{M}),\quad\forall x\in(-1,1)\,, \quad\lvert P''(x)\lvert\leq 4\Big(\frac m{\sqrt{1-x^{2}}}\Big)^{2}\lVert P\lVert_{\infty}\,.
\qe
\item[Laplace transform]  Let $0<\lambda_0<\lambda_1<\lambda_{2}<\dotsb$ be any real numbers. Let $\mathbf F_L=\left\{\exp(-\lambda_{0} x),\exp(-\lambda_{1} x),\exp(-\lambda_{2} x),\dotsc\right\}$. From Newman's inequality (\cite{borwein1995polynomials}, p.276), we know that:
\[\forall P\in\mathrm{Span}(\mathbf F_{L}),\quad\forall x\in[0,+\infty)\,, \quad\lvert P''(x)\lvert\leq \Big({9\sum_{i=0}^{\infty}\lambda_{i}}\Big)^{2}\lVert P\lVert_{\infty}\,.\]
\item[M\"untz polynomials I] Let $0<\alpha_1<\dotsb<\alpha_{m}$ be any real numbers. Let $\mathbf F_{M\ddot{u}}=\{1,x^ {\alpha_1},\dotsc,x^ {\alpha_m}\}$ then for all $\eta>0$, there is a constant $c_{\eta}$ such that
\[\forall P\in\mathrm{Span}(\mathbf F_{M\ddot{u}}),\quad\forall x\in(0,1-\eta)\,, \quad\lvert P''(x)\lvert\leq \frac {c_{\eta}}{x^{2}}\lVert P\lVert_{\infty}\,.\]
\item[M\"untz polynomials II] Let $1<\alpha_0<\dotsb<\alpha_{m}$ be any real numbers. Let $\mathbf F_{M\ddot{u}}=\{x^{\alpha_{0}},x^ {\alpha_1},\dotsc,x^ {\alpha_m}\}$ then, by Newman's inequality (\cite{borwein1995polynomials}, p.276), 
\[\forall P\in\mathrm{Span}(\mathbf F_{M\ddot{u}}),\quad\forall x\in(0,1]\,, \quad\lvert P''(x)\lvert\leq \frac {3^{4}(\sum_{i=0}^{m}\alpha_{i})(\sum_{i=0}^{m}\alpha_{i}+1)}{x^{2}}\lVert P\lVert_{\infty}\,.\]
\end{description}

\ni
These examples shows that the BIP property is a mild assumption in most practical cases.

\section{Support detection from noisy Fourier/moment samples}\label{sec:FourierCase}

\subsection{Detection from noisy Fourier samples}
In this subsection, we mention the example of Fourier samples to illustrate our results. Recently, much emphasis has been put on the recovery of a discrete measure from noisy band-limited data \cite{candes2012super,fernandez2013support,tang2013near}. In this setting, we observe noisy Fourier samples up until a frequency cut-off $f_c\in\mathbb N^\ast$. We specify notation:

\begin{itemize}
 \item The number of samples is $2f_c+1$ hence $m=2f_c$.
 \item For sake of simplicity, we place ourselves on $\mathbb T=[0,1]$.
 \item For all $k\in\{-f_c,\ldots,f_c\}$, we set for all $x\in[0,1]$, $\varphi_k(x)=\exp(\mathrm{i}\,2\pi kx)$, and $\Phi=(\varphi_{-f_c},\ldots,\varphi_{f_c})$.
 \item Assume $(\varepsilon_k)_{k=0}^{m}$ are random complex Gaussian:
 \[
  \varepsilon_k=\varepsilon_k^{(1)}+ \mathrm{i}\,\varepsilon_k^{(2)}\,,
 \]
where $\varepsilon_k^{(1)}$, $\varepsilon_k^{(2)}$, $k\in\{-f_c,\ldots,f_c\}$ are i.i.d. centered Gaussian random variables with standard deviation $\sigma$:
\[
 \varepsilon_k^{(1)}\sim \varepsilon_k^{(2)}\sim\mathcal{N}(0,\sigma^2)\,.
\]
\ni
We mention that $\varepsilon=(\varepsilon_{-f_c},\ldots,\varepsilon_{f_c})$. 
\item Finally, we recall that we observe $y=\int_{\mathbb T}\Phi\,\d\Delta+\varepsilon$.
\end{itemize}

\ni
Our results show that if the spikes are sufficiently separated, at least $2.5/f_c$ apart, then one can detect some point sources with a known precision solving a simple convex optimization program.

\begin{defn}[Minimum separation \cite{candes2012towards}] For a family of points
$\mathbf S\subset\mathbb T$, the minimum separation is defined as the closest distance
between any two elements from $\mathbf S$:
 \[
  \ell(\mathbf S)=\inf_{\substack{(x,x')\in \mathbf S^2\\ x\neq x'}}\lvert x-x'\lvert\,.
 \]
We emphasize that the distance is taken around the circle so that, for
example $\lvert 5/6-1/6\lvert=1/3$.
\end{defn}

\begin{figure}[h]
 \begin{center}
  \includegraphics[height=7cm,width=11.3262cm]{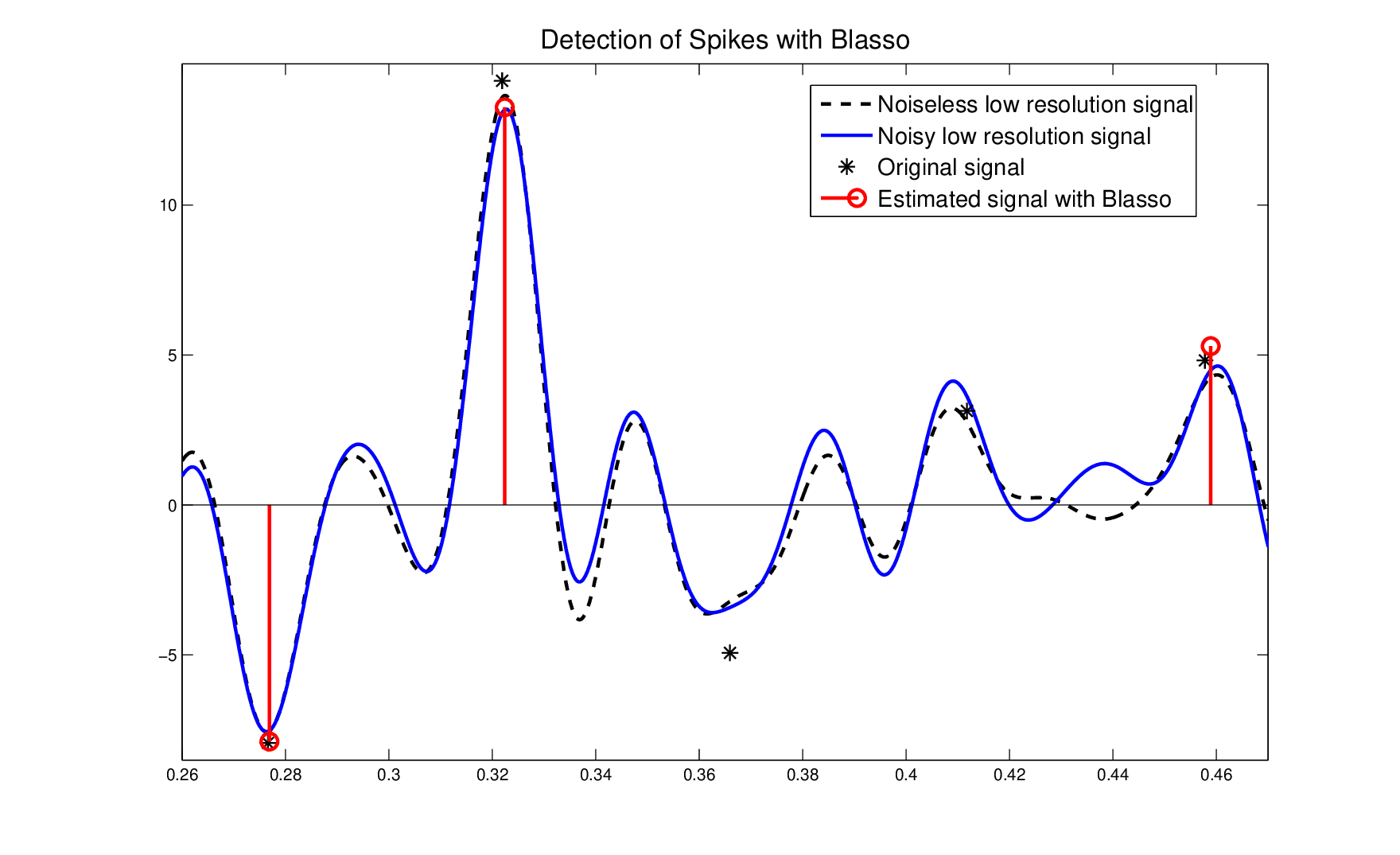}
 \end{center}
\caption{The problem is the following: we aim at recovering some spikes of the original signal (stars $\ast$) from the observation of a corrupted optical device (blue line) which can differ heavily from the true noiseless observation (black dotted line). Our procedure (red circles) provides a close estimate of the location of some spikes.}\label{Fig:Fourier}
\end{figure}

\ni
In this framework, we have the following theorem that quantifies the support detection error of BLASSO.

\begin{cor}\label{cor:Fourier} Assume that $f_{c}\geq128$. Let $\Delta$ be a discrete measure such
that:
\eq\label{eq:minimalseparationhypothesis}
 \ell(\mathbf S)\geq\dfrac{2.5}{f_c}\,,
\qe
where $\mathbf S$ denotes the support of $\Delta$. Let $\hat\Delta$ be a solution to
\eqref{def:Blasso} with tuning parameter $\lambda$ such that:
\[
 \lambda\geq\lambda_F:=2\,\sigma\,\sqrt{6\,f_c\,\log(f_c)}\,.
\]
then, with probability greater than $1-2\exp[-(\log f_c)\,{\lambda^2}/{\lambda_F^2}]$, the following holds. 
\begin{itemize}
\item For all $x\in [0,1]$ such that: 
\eq\label{eq:Fourierpoidsmin}
 \lvert\hat\Delta(\{x\})\lvert\geq 218\,\lambda\,,
\qe
there exists a unique $T\in\mathbf S$ satisfying:
\[
 \lvert x-T\lvert\leq \Big[\dfrac{\lambda}{0.1678\,\lvert\hat\Delta(\{x\})\lvert}\Big]^{1/2}\,\dfrac1{f_c}\leq\dfrac{0.1649}{f_c}\,,
\]

\item For all $T_{j}\in \mathbf S$ corresponding to a weight such that: 
\[
 \Delta_{j}>218\,\lambda\,,
\]
there exists $\hat T\in\mathrm{Supp}(\hat\Delta)$ satisfying:
\[
 \lvert T_{j}-\hat T\lvert\leq \Big[\dfrac{\lambda}{0.1678(\Delta_{j}-218\lambda)}\Big]^{1/2}\,\dfrac1{f_{c}}\,.
\]

\end{itemize}
\end{cor}
%
 \ni
This result provides a quantitative estimate of the location of spikes in terms of the output amplitudes or the input amplitudes. 
\begin{rem}
Observe that our procedure do not suppose any knowledge on the total number of spikes. This property is of great importance in actual practice. Only the minimal distance between any pair of atoms is relevant to BLASSO.
\end{rem}

\begin{rem}
The constant $2.5$ in the minimal separation assumption \eqref{eq:minimalseparationhypothesis} is not optimal. However, it can be lowered to $2$ using the result \cite{candes2012towards} but the constants in Corollary \ref{cor:Fourier} seems too large in this case. For instance, the constant $218$ in \eqref{eq:Fourierpoidsmin} would be greater than $10^{3}$. So, we choose to keep $2.5$ in \eqref{eq:minimalseparationhypothesis} so as to get not too large constants.
\end{rem}

\subsection{Detection from noisy moment samples}

In this section, we investigate the diffraction limit from moment measurements. In particular, we try to understand how BLASSO can super-resolve spikes from noisy moment samples. 

\begin{itemize}
 \item The number of samples is $m+1$.
 \item For sake of simplicity, consider that $\mathbb T=[-1,1]$.
 \item Consider $\mathbf F_{C}=\{\varphi_{0},\ldots,\varphi_{m}\}$ defined by $T_{0}=1$ and for all $k\in\{1,\ldots,m\}$,
\eq\label{eq:ChebyshevFamily}
\varphi_{k}=\sqrt 2\,T_{k}\,,
\qe
where $T_{k}(x)=\cos(k\arccos(x))$ is the $k$-th Chebyshev polynomial of the first kind. Hence, the family $\mathbf F_{C}$ is an orthonormal family with respect to the probability measure $\Pi(\d t)=(1/\pi)\,({1-t^{2}})^{-1/2}\,\mathcal L(\d t)$ on $[-1,1]$ where $\mathcal L$ denotes the Lebesgue measure.
 
  \item Assume $\varepsilon_k$ are i.i.d. centered Gaussian random variables with standard deviation $\sigma$. We mention that $\varepsilon=(\varepsilon_{0},\varepsilon_{1},\ldots,\varepsilon_{m})$. 
\item Finally, we recall that we observe $y=\int_{\mathbb T}\Phi\,\d\Delta+\varepsilon$.
\end{itemize}
\begin{cor}\label{cor:Moment} Let $m\geq 9$. Let $\Delta$ be a discrete measure such that its support $\mathbf S$ is $2c_{0}$ far from the endpoints of $\T$, namely $\pm1$. Assume that $\Delta$ satisfies $\mathrm{QIC}(C_{a},C_{b})$ with ${C_{b}/C_{a}}\leq c_{0}^{2}$. Let $\hat\Delta$ be a solution to
\eqref{def:Blasso} with tuning parameter $\lambda$ such that:
\[
 \lambda\geq\lambda_M:=\sigma\sqrt{6m\log m}\,.
\]
then, with probability greater than $1-8\frac\lambda{\lambda_{M}}\sqrt{\log m}\,\exp\big[-(2\frac{\lambda^{2}}{\lambda_M^{2}}-1)\log m\big]$, the following holds. 
\begin{itemize}
\item For all $x\in [-1,1]$ such that: 
\[
 \lvert\hat\Delta(\{x\})\lvert>\frac{2\lambda}{C_{b}}\,,
\]
there exists a unique $T\in\mathbf S$ satisfying:
\[
 \lvert x-T\lvert\leq \Big[\dfrac{2\lambda}{C_{a}\,\lvert\hat\Delta(\{x\})\lvert}\Big]^{1/2}\,\dfrac1{m}\leq\dfrac{c_{0}}{m}\,,
\]

\item Let $C'=2+\max\Big\{\frac{2(1-C_{b})}{C_{b}},\frac{4}{C_{a}-C_{b}}\Big\}$. Then, for all $T_{j}\in \mathbf S$ corresponding to a weight such that: 
\[
 \Delta_{j}>C'\,\lambda\,,
\]
there exists $\hat T\in\mathrm{Supp}(\hat\Delta)$ satisfying:
\[
 \lvert T_{j}-\hat T\lvert\leq \Big[\dfrac{2\lambda}{C_{a}(\Delta_{j}-C'\lambda)}\Big]^{1/2}\,\dfrac1{m}\,.
\]

\end{itemize}
\end{cor}

\section{Rice method}
\ni
We use the Rice method as a convenient tool to upper bound the supremum of a random polynomial on a compact set. Proofs can be found in \ref{sec:ProofRice}.

\subsection{Polynomial case} Consider the family $\mathbf F_{C}=\{\varphi_{0},\ldots,\varphi_{m}\}$ defined by \ref{eq:ChebyshevFamily}. Define the Gaussian process $\{X_m(t),\ t \in[-1,1]\}$ by:
\[
  \forall t\in[-1,1],\quad X_m(t) =  \xi_0\varphi_{0}(t) + \xi_1 \varphi_{1}(t) +  \xi_2 \varphi_{2}(t) +  \ldots +  \xi_m \varphi_{m}(t)\,,
\]
 where $ \xi_1,  \ldots ,  \xi_m$ are i.i.d. standard normal. Its covariance function is 
 \[
 r(s,t) = 1 + \varphi_{1}(t)\varphi_{1}(s) + \varphi_{2}(t)\varphi_{2}(s) + \ldots +  \varphi_{1}(m)\varphi_{m}(s)\,,
 \]
where the dependence in $m$ has been omitted. Observe that its maximal variance is attained at point $1$ and is given by $\sigma^2_m = 2m+1$, and its variance function is $\sigma^2_m (t) =  1 + \varphi_{1}(t)^{2} + \varphi_{2}(t)^{2} + \cdots +  \varphi_{1}(m)^{2}$.
 
\begin{prop}\label{prop:RicePoly}
Let $\displaystyle M = \max_{t \in [-1,1]} | X_m(t)|$. Then, for $m\geq 12$ and for $u > \sqrt{1+2m}$,
\[
  \P\{ M>u\}  \leq   \frac{4 m (1+u)}{\sqrt{2\pi}}  \exp\big(-\frac{u^{2}}{1+2m}\big)\,.
\]
\end{prop}

 \ni

\subsection{Fourier case} We consider  the trigonometric functions:
\[
\varphi_k (t)= \exp(  i 2 \pi k t  )\,,\quad t \in [0,1] \ \mathrm{and}\  k  \in\mathcal{ K}  := \{ -f_c, \ldots,f_c\}\,,
\]
and random  complex Gaussian errors:
    $$
     \varepsilon_k = \varepsilon^{(1)}_k + \mathrm{i}\, \varepsilon^{(2)}_k , 
     $$
      where  the variables  $\varepsilon^{(1)}_k , \varepsilon^{(2)}_k,  k   \in\mathcal{ K} $ are independent with standard  normal distribution.
       
\begin {prop}\label{prop:RiceFourier}
Let $\displaystyle Z(t)   =  \sum_{k   \in\mathcal{ K} }  \varepsilon_k    \varphi_k (t)$. Then, for $ u>\sqrt{2} $,
      \[
       \P\{ \sup_{t \in [0,1]}  \| Z(t) \| >u \}\leq  4 \big( 
        \exp( -\frac{u^2}{2(2{f_c}+1)})  +  \sqrt{ \frac{{f_c}  ({f_c}+1)}{3}}  \exp\big ( - \frac{u^2}{4(2{f_c}+1)}\big)\big)\,.
       \]
      \end{prop}
      \ni

\appendix
\section{Proof of Theorem \ref{thm:LocQICGeneral}}\label{sec:ProofQIC1}
\ni
Let $\displaystyle P=\sum_{k=0}^m a_k \varphi_k$ given by the QIC condition such that $P$ is a dual certificate of $\Delta$. Set: \[d=\lVert{\hat\Delta}\lVert_{TV}-\lVert{\Delta}\lVert_{TV}-\Re(\int_{\mathbb T} P\,\d(\hat\Delta-\Delta))\,,\]
where $\Re(\,.\,)$ denotes the real part. One can check that $d$ belongs to the Bregman divergence of the TV-norm between $\hat\Delta$ and $\Delta$. Hence, we know that it is non-negative. Indeed, observe that:
\[
\Re(\int_{\mathbb T} P\,\d(\hat\Delta-\Delta))=\Re(\int_{\mathbb T} P\,\d(\hat\Delta))-\lVert{\Delta}\lVert_{TV}\leq\lVert{\hat\Delta}\lVert_{TV}-\lVert{\Delta}\lVert_{TV}
\]
since $P$ is sub-gradient of the TV-norm at point $\Delta$. From the definition of BLASSO, we have:
\[
 \frac12\lVert\int_{\mathbb T}\Phi\,\d\hat\Delta-y\lVert_2^2+\lambda\lVert\hat\Delta\lVert_{TV}\leq\frac12\lVert\varepsilon\lVert_2^2+\lambda\lVert \Delta\lVert_{TV}\,.
\]
Therefore:
\eq\label{eq:resultintermediaire}
 \frac12\lVert\int_{\mathbb T}\Phi\,\d\hat\Delta-y\lVert_2^2+\lambda\,d+\lambda\,\Re(\int_{\mathbb T} P\,\d(\hat\Delta-\Delta))\leq\frac12\lVert \varepsilon\lVert_2^2\,.
\qe
It follows that:
\[
 \frac12\lVert\int_{\mathbb T}\Phi\,\d\hat\Delta-y+\lambda\,\overline{a}\lVert_2^2+\lambda\,d\leq\frac12\lVert\lambda\,{a}\lVert_2^2+\frac12\lVert \varepsilon\lVert_2^2-\lambda\,\langle\varepsilon\,,\, \overline{a }\rangle\,.
\]
Eventually, we get:
\[
 d\leq\frac\lambda2\big\lVert\overline{a}-\frac{{\varepsilon}}\lambda\big\lVert_2^2\,.
\]
Using the triangular inequality and Parseval's identity (we recall that $\mathbf F$ is an orthonormal family with respect to a probability measure $\Pi$ on $\mathbb T$), it yields:
\begin{align*}
  \big\lVert\overline a-\frac{{\varepsilon}}\lambda\big\lVert_2&\leq\big\lVert a\lVert_2+\lVert\frac{{\varepsilon}}\lambda\big\lVert_2\,,\\
  & =\Big(\int_{\mathbb T}\lvert P\lvert^2(x)\Pi(\d x)\Big)^{1/2}+\dfrac{\lVert\varepsilon\lVert_{2}}\lambda\,.
\end{align*}
Since $\lambda\geq\lVert\varepsilon\lVert_{2}$ and $\lVert P\lVert_\infty\leq1$, we have:
\eq\label{eq:majfourier}
 \big\lVert\overline a-\frac{{\varepsilon}}\lambda\big\lVert_2\leq2\,.
\qe
\ni
So we get that $d\leq 2\lambda$. Moreover, observe that:
\begin{align}
 d&=\lVert\hat\Delta\lVert_{TV}-\Re(\int_{\mathbb T}P\ \d\hat\Delta)\,,\notag
 \\&=\sum_{k=1}^n\hat\Delta_k\,\big[1-\lvert P\lvert(\hat t_k)\,\cos(\hat\theta_k+\theta_P(\hat T_k))\big]\,,\label{eq:BregmannPhase}\\
 & \geq\sum_{k=1}^n\hat\Delta_k\min\{C_{a}m^{2}\min_{T\in\mathbf S}d(\hat T_{k},T)^{2},C_{b}\}\,.\label{eq:WeakWasserstein}
\end{align}
This shows the results $(1)$ and $(2)$. Observe that the localization estimate in the statement of the theorem is a consequence of the aforementioned inequality.


\section{Proof of Theorem \ref{thm:LocBIPGeneral}}
\ni
We begin with a lemma on the optimality conditions of Blasso. For all $a\in\mathbb C^{m+1}$ and for all $\mu\in\mathcal M$, we denote 
\[
\langle a\,,\,\Phi\rangle=\sum_{k=0}^{m}\overline a_{k}\varphi_{k}\quad\mathrm{and}\quad c(\mu)=\int_{\mathbb T}\Phi\,\d\mu,.
\]
\ni
We recall some standard lemmas. The proofs are given in Appendix \ref{sec:AppendixAuxLem}.
\begin{lem}\label{lem:Optimality}
Any solution $\hat\Delta$ to \eqref{def:Blasso} satisfies the following optimatlity conditions:
\begin{enumerate}
\item
$\lVert \langle c(\hat\Delta)-y\,,\,\Phi\rangle\lVert_{\infty}\leq\lambda$;
\item
$\langle y-c(\hat\Delta)\,,\,c(\hat\Delta)\rangle=\lambda\lVert \hat\Delta\lVert_{TV}$.
\end{enumerate}
\end{lem}
\begin{lem}\label{lem:MajorationL1}
Set $\nu=\hat\Delta-\Delta$ then the following inequality holds:
\eq\label{eq:PredicitonL1}
\forall P\in\mathrm{Span}(\mathbf F),\quad\lvert\int_{\mathbb T}P\,\d\nu\lvert\leq(\lambda+\lambda_{0})\lVert P\lVert_{1}\,,
\qe
where $\lambda_{0}=\lVert \langle \varepsilon\,,\,\Phi\rangle\lVert_{\infty}$.
\end{lem}
\ni 
Using the BIP property and the QIC condition, we prove the following interpolation lemma.
\begin{lem}\label{lem:InterpolationOnePhase}
Let $j\in\{1,\ldots,s\}$ then there exists a generalized polynomial $Q_{j}\in\mathbf F$ satisfying the following properties:
 \begin{itemize}
  \item $\forall k\in\{1,\ldots,s\}\,,\quad Q_{j}(T_k)=\delta_{k,j}$,
  \item $\forall x\neq T_{j}\,,\quad \lvert Q_{j}(x)\lvert<1$,
  \item $\forall x\in\T$ such that $d(x,T_{j})\leq c_{0}/m$, it holds:
\[\lvert 1- Q_{j}(x)\lvert\leq (C_{c}/2)\,m^2\,d(x,T_{j})^2\,,\]
 \item $\forall k\in\{1,\ldots,s\}\setminus\{j\}$, $\forall x\in\T$ such that $d(x,T_{k})\leq c_{0}/m$, it holds:
 \[\lvert Q_{j}(x)\lvert\leq (C_{c}/2)\,m^2\,d(x,T_{k})^2\,,\]
\item $\forall x\in\mathcal F(c_{0},\mathbf S)\,,\quad\lvert Q_{j}(x)\lvert\leq1-C_b$.
 \end{itemize}
\end{lem}
\ni
Using these interpolating polynomials and \eqref{eq:WeakWasserstein} we get that for all $j\in\{1,\ldots,s\}$:
\begin{align*}
&\lvert \sum_{\{k:\ d(\hat T_{k},T_{j})>\frac{c_{0}}m\}}\hat\Delta_{k}\,\exp(\mathrm{i}\,\hat\theta_k)\,Q_{j}(\hat T_{k})+\sum_{\{k:\ d(\hat T_{k},T_{j})\leq\frac{c_{0}}m\}}\hat\Delta_{k}\,\exp(\mathrm{i}\,\hat\theta_k)\,(Q_{j}(\hat T_{k}) -1)\lvert \,,\\
&\leq  \sum_{\{k:\ d(\hat T_{k},T_{j})>\frac{c_{0}}m\}}\hat\Delta_{k}\,\lvert Q_{j}(\hat T_{k})\lvert+\sum_{\{k:\ d(\hat T_{k},T_{j})\leq\frac{c_{0}}m\}}\hat\Delta_{k}\,\lvert Q_{j}(\hat T_{k}) -1\lvert \,,\\
&\leq \sum_{k=1}^n\hat\Delta_k\min\{(C_{c}/2)m^{2}\min_{T\in\mathbf S}d(\hat T_{k},T)^{2},1-C_{b}\}\,,\\
& \leq \max\Big\{\frac{1-C_{b}}{C_{b}},\frac{C_{c}}{2C_{a}}\Big\}\times\sum_{k=1}^n\hat\Delta_k\min\{C_{a}m^{2}\min_{T\in\mathbf S}d(\hat T_{k},T)^{2},C_{b}\}\,,\\
& \leq \lambda \max\Big\{\frac{2(1-C_{b})}{C_{b}},\frac{C_{c}}{C_{a}}\Big\}\,.
\end{align*}
\ni
Now, invoking \eqref{eq:PredicitonL1} we get:
\begin{align*}
\lvert \Delta_{j}-\sum_{\{k:\ d(\hat T_{k},T_{j})\leq\frac{c_{0}}m\}} \hat\Delta_{k}\lvert&=\lvert \int_{\T}Q_{j}\,\d\Delta-\int_{\T}Q_{j}\,\d\hat\Delta+\sum_{\{k:\ d(\hat T_{k},T_{j})>\frac{c_{0}}m\}}\hat\Delta_{k}\,\exp(\mathrm{i}\,\hat\theta_k)\,Q_{j}(\hat T_{k})\\
&+\sum_{\{k:\ d(\hat T_{k},T_{j})\leq\frac{c_{0}}m\}}\hat\Delta_{k}\,\exp(\mathrm{i}\,\hat\theta_k)\,(Q_{j}(\hat T_{k}) -1)\lvert \,,\\
& \leq (2+\max\Big\{\frac{2(1-C_{b})}{C_{b}},\frac{C_{c}}{C_{a}}\Big\})\lambda\,.
\end{align*}
\ni
This proves $(3)$. From $(1)$ and $(2)$ given by Theorem \ref{thm:LocQICGeneral}, we deduce the last part of the theorem.

\section{Proofs of the auxiliary lemmas}\label{sec:AppendixAuxLem}

   \subsection{Proof of Lemma \ref{lem:Optimality}}
The convex function $f(\nu)=\frac12\lVert
\int_{\mathbb T}\Phi\,\d\nu-y\lVert^2_2+\lambda\lVert\nu\lVert_{TV}$ is minimized at point $\hat\Delta$ then $0\in\partial f(\hat\Delta)$ or equivalently, for all $\nu\in\mathcal M$,
\begin{equation}\label{eq:optimality2}
\lambda(\lVert \nu\lVert_{TV}-\lVert \hat\Delta\lVert_{TV})\geq\,\langle y-c(\hat\Delta)\,,\,c(\nu-\hat\Delta)\rangle\,.
\end{equation}
Conversely, we can easily check that if \eqref{eq:optimality2} holds then $\hat\Delta$ is a minimizer. Thus, \eqref{eq:optimality2} is necessary and sufficient for $\hat\Delta$ to minimize $f(\nu)$. Observe that \eqref{eq:optimality2} leads to:
\begin{equation}\label{eq:optimality3}
\lambda\lVert \hat\Delta\lVert_{TV}-\langle y-c(\hat\Delta)\,,c(\hat\Delta)\rangle\leq\inf_{\nu\in\mathcal M}\{\lambda\lVert \nu\lVert_{TV}-,\langle y-c(\hat\Delta)\,,c(\nu)\rangle\}\,.
\end{equation}
Now, observe that, for all $z\in\mathbb C^{m+1}$,
\[
\sup_{\nu\in\mathcal M}\{\langle z\,,c(\nu)\rangle-\lVert \nu\lVert_{TV}\}=\sup_{\nu\in\mathcal M}\{\int_{\T}\langle  z\,,\Phi\rangle\d\nu-\lVert \nu\lVert_{TV}\}=\begin{cases}
0 &  \mathrm{if\ }\lVert\langle z\,,\,\Phi\rangle \lVert_{\infty}\leq 1;\\
\infty & \text{otherwise.}
\end{cases}
\]
Using this in \eqref{eq:optimality3}, we find that $\hat\Delta$ is a minimizer if and only if $(1)$ and $(2)$ hold.
\subsection{Proof of Lemma \ref{lem:MajorationL1}}
Let $\overline a$ be the coefficients of $P$. It holds
\begin{align*}
\int_{\T}P\,\d\nu&=\int_{\T}\langle a\,,\,\Phi\rangle\,\d\nu\,,\\
&=\langle a\,,\,\int_{\T}\Phi\,\d\nu\rangle_{\ell^{2}(\mathrm{C}^{m+1})}\,,
\\
&=\langle a\,,\,c(\hat\Delta)-c(\Delta)\rangle_{\ell^{2}(\mathrm{C}^{m+1})}\,,
\\
&=\langle P\,,\,E\rangle_{L^{2}_{\Pi}(\T)}\,,
\end{align*}
using the Parseval identity where $E$ denotes the trigonometric polynomial of order $n$ such that $E=\langle c(\hat\Delta)-x(\Delta)\,,\,\Phi\rangle$. By H\"older's inequality we have
\begin{equation*}
\lvert\int_{\T}P\,\d\nu\lvert\leq \lVert P\lVert_{1}\lVert E\lVert_{\infty}\,.
\end{equation*}
Now, we have to upper bound $\lVert E\lVert_{\infty}$. By triangular inequality, we have
\begin{equation*}
\lVert E\lVert_{\infty}\leq\lVert \langle c(\hat\Delta)-y\,,\,\Phi\rangle\lVert_{\infty}+\lVert \langle \varepsilon\,,\,\Phi\rangle\lVert_{\infty} \,.
\end{equation*}
By Lemma \ref{lem:Optimality}, we know that the first term in the right hand side is upper bounded by $\lambda$. By hypothesis, we know that the second term in the right hand side is upper bounded by $\lambda_{0}$.
\subsection{Proof of Lemma \ref{lem:InterpolationOnePhase}}
Without loss of generality, assume that $j=1$ and set $\theta_{1}=0$, $\tilde\theta_{1}=1$, $\theta_{k}=0$ and $\tilde\theta_{k}=\pi$, for all $k\in\{2,\ldots,s\}$. Invoke QIC to get $P$ and $\tilde P$ such that:
\begin{itemize}
  \item $P(T_1)=\tilde P(T_{1})=1$,
  \item $\forall k\in\{2,\ldots,s\}\,,\quad P(T_k)=-\tilde P(T_{k})=1$,
\item $\forall x\in\mathcal F(c_{0},\mathbf S)\,,\quad\max(\lvert P(x)\lvert,\lvert\tilde P(x)\lvert)\leq1-C_b$.
 \end{itemize}
 Set $Q_{1}=(1/2)(P+\tilde P)$ then one can check that:
  \begin{itemize}
  \item $\forall k\in\{1,\ldots,s\}\,,\quad Q_{1}(T_k)=\delta_{k,1}$,
  \item $\forall t\neq T_{1}\,,\quad \lvert Q_{1}(t)\lvert<1$,
  \item $\forall x\in\mathcal F(c_{0},\mathbf S)\,,\quad \lvert Q_{1}(x)\lvert\leq1-C_b$,
 \end{itemize}
 where $\delta_{k,1}$ equals $0$ if $k\neq1$ and $1$ if $k=1$. Moreover, since $P$ and $\tilde P$ reach their maxima at each point of $T$, their first derivative vanishes on $T$. Therefore, $Q_{1}'$ vanishes on $T$. Using a Taylor's theorem with explicit remainder and BIP, we conclude the proof.

\section{Proof of Corollary \ref{cor:Fourier}}\label{sec:proofthmfourier}

\ni
We begin with a key result.
\begin{thm}[Lemma 2.5 in \cite{candes2012towards}]\label{lem:Candes}
Let $\mathbf S=\{T_1,\ldots,T_s\}\subset[0,1]$ be the support of the target measure
$\Delta$. If $\ell(\mathbf S)\geq 2.5/{f_c}$ then there exists $P\in\mathrm{Span}(\mathbf F)$ such that for all $k=1,\ldots,s$:
 \begin{itemize}
  \item $P(T_k)=\exp(-\mathbf{i}\theta_k)$,
  \item Bound on the Taylor expansion at point $T_k$:
\[\forall x\in\Big[T_k-\dfrac{0.3313}{m},\,T_k+\dfrac{0.3313}{m}\Big]\,,
\quad\lvert
P(x)\lvert\leq1-0.0839\,m^2\,(x-T_k)^2\,,\]
\item Bound on the complement:
 \[\forall
x\in[0,1]\setminus\bigcup_{k=1}^s\Big[T_k-\dfrac{0.3313}{m},\,T_k+\dfrac{0.3313}{
m}\Big]\,,\quad\lvert P(x)\lvert\leq1-0.0092\,.\]
 \end{itemize}
\end{thm}
\ni
Hence, the support $\mathbf S$ enjoys $\mathrm {QIC}(0.0838,0.0092)$ and Theorem \ref{thm:LocQICGeneral} applies. Moreover, the support $\mathbf S$ satisfies $\mathrm{BIP}(1,\pi^{2})$, see \eqref{eq:BernsteinFourier}, and Theorem \ref{thm:LocBIPGeneral} can be invoked. Using the identity, for $a>0$ and $b>0$,
\eq
\forall u^{2}\geq2ab\,,\quad -\frac{u^{2}}{a}+b\leq-\frac{u^{2}}{2a}\,,
\qe
and Propositon \ref{prop:RiceFourier}, we get that $\lambda\geq\lVert \sum_{k=-f_c}^{f_c} \overline{\varepsilon_k}\varphi_k\lVert_\infty$ with a probability described in the statement of Corollary \ref{cor:Fourier}. 

\section{Proof of Corollary \ref{cor:Moment}}\label{sec:proofthmmoment}
\ni
Observe that $\mathbf S$ satisfies $\mathrm{BIP}(c_{0},4/(1-c_{0}^{2}))$ thanks to \eqref{eq:BernsteinMoment}. Take 
\[
u=\sqrt{6m\log m}\,\lambda/\lambda_{0}\,,
\]
in Propositon \ref{prop:RicePoly}. Hence, we get that $\lambda\geq\lVert \sum_{k=0}^{m} \overline{\varepsilon_k}\varphi_k\lVert_\infty$ with a probability described as in the statement of Corollary \ref{cor:Moment}. Therefore, Theorems \ref{thm:LocQICGeneral} and \ref{thm:LocBIPGeneral} give the result.

\section{Rice formula}\label{sec:ProofRice}
\subsection{Polynomial case}
By the change of variables $t=\cos\theta$, for all $t\in[-1,1]$:
\[
X_m(t) =X_{m}(\cos\theta)=  \xi_0 + \sqrt2\,\xi_1 \cos(\theta) +  \sqrt2\,\xi_2 \cos(2\theta) +  \ldots +  \sqrt2\,\xi_m \cos(m\theta) \,.
\]
Set $T_{m}(\theta) :=X_{m}(t)$. We recall that its variance function is given by:
\[\sigma^{2}_m(\theta)=1+2\,\cos^{2}(\theta) + 2\,\cos^{2}(2\theta) +  \ldots +  2\,\cos^{2}(m\theta)=\frac{1+2m+\mathbf D_{m}(2\theta)}2\,,\]
where $\mathbf D_{m}$ denotes the Dirichlet kernel of order $m$. Observe that:
 \[\forall \theta\in\R\,,\quad\frac{1+2m}3=\sigma^{2}_m(3\pi/(2+4m))\leq \sigma^{2}_m(\theta)\leq 1+2m\,,\]
so that $T_m\big(\frac{3\pi}{2(1+2m)}\big)\sim\mathcal N(0,\frac{1+2m}3)$. By the Rice method \cite{Aza20081190}, for $u>0$:
 \begin{align*}
 \P\{M>u\}& \leq 2  \P\{   \max_{\theta \in [0,\pi]} T_m(\theta)  >u\} \,,\\
& \leq
 2  \P\{T_m(0)>u\} +2\, \E[ U_u([0,\pi])]\,,\\
 &=   2[1-\Psi\Big(\frac{u}{\sqrt{1+2m}}\Big)]+  2 \int_{0}^\pi \E\big( (T'_m(\theta))^+ \big| T_m(\theta) =u) \psi_{\sigma_m(\theta)} (u)\d t\,,
 \end{align*}
 where $U_u$ is the number of crossings of the level $u$, $\Psi $ is the standard normal distribution, and $\psi_{\sigma}$ is the density of the centered normal distribution with standard error $\sigma$. First, observe that for $v>0$, $(1-{\Psi}(v) )\leq \exp(-v^2/2) $. Hence, 
 \[
 1-\Psi\Big(\frac{u}{\sqrt{1+2m}}\Big)\leq\exp\big(-\frac{u^{2}}{1+2m}\big)\,.
 \]
 Moreover, regression formulas implies that:
 \begin{align*}
  \E\big( T'_m(\theta)\big| T_m(\theta) =u\big) & =\frac{ r_{0,1} (\theta,\theta)}{r(\theta,\theta)} u\,, \\
  \Var\big( T'_m(\theta)\big| T_m(\theta) =u\big) & \leq   \Var\big( T'_m(\theta)\big) =  r_{1,1} (\theta,\theta)\,,
  \end{align*}
  where, for instance,  $r_{1,1} (\nu,\theta) =\frac{ \partial^2 r(\nu,\theta)}{\partial \nu \partial \theta}$. We recall that the covariance function is given by:
  \begin{align*}
   r(\nu,\theta)&=1+2\,\cos(\nu)\cos(\theta) + 2\,\cos(2\nu)\cos(2\theta) +  \ldots +  2\,\cos(m\nu)\cos(m\theta)\,, \\
   & = \frac 12 \big[\mathbf D_{m}(\nu-\theta)+\mathbf D_{m}(\nu+\theta) \big]\,.
  \end{align*}
Observe that:
 \begin{align*}
  r_{0,1} (\theta,\theta) & = \frac12\,\mathbf D'_{m}(2\theta)=-\sum_{k=1}^{m}k\sin(2k\theta)\,,\\
   \mathrm{and}\quad r_{1,1} (\theta,\theta) & = \frac{(2m+1)(m+1)m}6-\sum_{k=1}^{m}k^{2}\cos(2k\theta)\,.
   \end{align*}
 On the other hand, if $Z \sim \mathcal N(\mu, \sigma^2)$ then
\[
 \E(Z^+ ) = \mu\, \Psi\big( \frac\mu\sigma\big) + \sigma\,  \psi\big( \frac\mu\sigma\big) \leq  \mu^+ + \dfrac\sigma{\sqrt{2 \pi}}\,,
\]
where $\psi$ is the standard normal density. We get  that:
 \begin{align*}
 &\int_{0}^\pi \E\big( (T'_m(\theta))^+ \big| T_m(\theta) =u) \psi_{\sigma_m(\theta)} (u)\d t\leq\int_{0}^\pi\frac{ [\mathbf D'_{m}(2\theta)]^{+}}{2\,\sigma^{2}_{m}(\theta)}\, u \psi_{\sigma_m(\theta)}(u) \d\theta\\&  + \frac{1}{\sqrt{2\pi}} \int_{0}^\pi\big[\frac{(2m+1)(m+1)m}6-\sum_{k=1}^{m}k^{2}\cos(2k\theta)\big]^{1/2} \,\psi_{\sigma_m(\theta)}(u) \d\theta\,, \\
  &= A+B\,.
\end{align*}
We use the following straightforward relations: 
\begin{itemize}
  \item $\forall\,0<\sigma_1 <\sigma_2 <u\,,\quad\psi_{\sigma_1}(u) \leq \psi_{\sigma_2}(u)$, 
  \item $\forall\theta\,,\quad[\mathbf D'_{m}(2\theta)]^{+}\leq \frac{m(m+1)}2$, 
    \item $\forall \theta\,,\quad\frac{1+2m}3\leq \sigma^{2}_m(\theta)$,
  \end{itemize}
Eventually, we get, for $u> \sqrt{2m+1}$:
\begin{align*}
A &\leq  \frac{3\pi}2\, m\,u  \psi_{\sqrt{2m+1}} (u)\,,\\
B&  \leq \Big[\frac\pi{12}(2m+1)(m+1)m\Big]^{1/2}\,\psi_{\sqrt{2m+1}} (u)\,,
\end {align*}
and therefore:
\[
\P\{ M>u\} \leq 2\exp\big(-\frac{u^{2}}{1+2m}\big)+2\Big[\frac{3\pi\,mu}{2\sqrt{2m+1}}+(\frac\pi{12}m(m+1))^{1/2}\Big]\,\psi\big(\frac u{\sqrt{2m+1}}\big)\,.
\]
Morover, for $m\geq 12$, $\pi m(m+1)/12\leq m^2$ and ${3\pi\,mu}/(2\sqrt{2m+1}) \leq mu$, yielding:
\[
\P\{ M>u\} \leq \frac{4 m (1+u)}{\sqrt{2\pi}}  \exp\big(-\frac{u^{2}}{1+2m}\big)\,.
\]
%
\subsection{Fourier case}
\ni
       We have:
       \begin{align*}
        Z(t) &= \varepsilon^{(1)} _0  + \sum_{k=1}^{f_c}  \big(  \varepsilon_k^{(1)} +  \varepsilon_{-k}^{(1)} \big)  \cos(2\pi k t)  + \big(  \varepsilon_k^{(2)} - \varepsilon_{-k}^{(2)} \big)  \sin(2\pi k t) 
        \\
        & + \mathrm{i}\, \Big[ 
      \varepsilon^{(2)} _0  +\sum_{k=1}^{f_c}  \big(  \varepsilon_k^{(2)} +  \varepsilon_{-k}^{(2)} \big)  \cos(2\pi k t)  + \big(  \varepsilon_k^{(1)} - \varepsilon_{-k}^{(1)} \big)  \sin(2\pi k t) \Big].
\end{align*}
One can see  that  $ Z(t) = X(t) + \mathrm{i}\, Y(t)$ where $X(t)$ and $Y(t)$ are two independent Gaussian stationary processes with the same auto-covariance function:
\[
 \Gamma(t) =  1 + 2 \sum_{k=1}^{f_c}    \cos(2\pi k t) = \mathbf D_{f_c}(t)\,,
\]
where $\mathbf D_{f_c}(t) $ denotes the Dirichlet kernel. Set:
\[
 {\sigma_m}^2  = \Var(X(t)) = \mathbf D_{f_c}(0) = 2{f_c} +1\,.
\]
We use the following inequalities:
\begin{align*}
 \P \{ \|Z\|_{\infty} >u\}  &\leq  \P \{  \|X\|_{\infty} >u/\sqrt{2}  \}+ \P\{ \|Y\|_{\infty} >u/\sqrt{2}\}\,, \\&=  2 \P \{  \|X\|_{\infty} >u/\sqrt{2}\} \,.
\end{align*}
and 
\begin{equation}\label{f:mm}
\P \{  \|X\|_{\infty} >u/\sqrt{2}   \}    \leq 2  \P \{      \sup_{t \in[0,1]}  X(t) >u/\sqrt{2} \}\,.
\end{equation}
To give bounds to the right hand side of (\ref{f:mm}), we use the Rice method \cite{Aza20081190} using the fact that the process $X(t)$ (for example)  is periodic with $\Gamma( \frac{1}{2{f_c}+1}) =0$:
\begin{align*}
 \P \{ \sup_{t \in[0,1] }X(t) >u/\sqrt{2} \}
&= \P \{ \forall  t \in[0,1] ;X(t) >u/\sqrt{2}\} + \P\{ U_{u/\sqrt{2} }>0\} \,,\\ &\leq  \bigg( \overline{\Psi} (u/ (\sqrt{2} {\sigma_m})) \bigg)^2 + \E (U_{u/\sqrt{2} })\,,
\end{align*}
where $U_v$ is the number of up-crossings of the level $v$ by the process $X(t)$ on the interval $[0,1]$ and $\overline{\Psi} $ is the tail of the standard normal distribution. By the Rice formula:
\[
  \E (U_{u/\sqrt{2} })  = \frac{1} {2\pi}          \sqrt{ \Var(X'(t)) }    \frac{1}{ {\sigma_m}}     \exp( -\frac   {u^2}   {4 {\sigma_m}^2}  )\,,
\]
where:
\[
 \Var(X'(t))   = -\Gamma''(0)  = 2 (2\pi)^2 \sum_{k=1}^{f_c} k^2 = \frac{4 \pi^2 }{3} {f_c} ({f_c}+1)(2{f_c}+1)\,.
\]
The following inequality is well known : for $v>0$, $\overline{\Psi}(v) \leq \exp(-v^2/2) $, it yields:
  \begin{equation*}\label{f:x}
 \P \{ \sup_{t \in[0,1] }X(t) >\dfrac u{\sqrt{2}} \}   \leq  \exp( -\frac{u^2}{2(2{f_c}+1)})  +  \sqrt{ \frac{{f_c}  ({f_c}+1)}{3}}  \exp\big ( - \frac{u^2}{4(2{f_c}+1)}\big)\,.
   \end{equation*}
   The result follows.

\subsection*{Acknowledgments}
The authors would like to thank anonymous referee for her/his fruitful comments and suggestions.

 \bibliographystyle{plain}
 \bibliography{biblio}
\end{document}